\documentclass{article}

\usepackage[utf8]{inputenc}
\usepackage[english]{babel}
\usepackage{amsmath,amsfonts,amssymb,amsthm}
\usepackage{mathrsfs}
\usepackage{graphicx}
\usepackage{hyperref}
\usepackage{url}

\DeclareMathOperator{\diam}{diam}
\newcommand{\R}{\mathbb{R}}

\newtheorem{theorem}{Theorem}

\newtheorem{claim}{Claim}
\newtheorem{question}{Question}
\newtheorem{definition}{Definition}
\newtheorem{conjecture}{Conjecture}
\newtheorem{problem}{Problem}

\title{A counterexample to Gale's conjecture in $\mathbb{R}^3$}
\author{Vsevolod Voronov}
\date{20 October 2026}

\begin{document}

\maketitle

\begin{abstract}
It is known that for every partition of the sphere $S^2$ of unit diameter into four parts, the largest diameter of a part is at least $\sqrt{(3+\sqrt{3})/6} = 0.88807\dots$, and this bound is attained. The conjecture of D.~Gale (1953) asserts that no set in $\mathbb{R}^3$ has a larger ratio. In this paper, we construct a $101$-point set $V\subset\mathbb{R}^3$ such that every partition of $V$ into four parts contains a part of diameter greater than $0.88928\,\diam V$.
\end{abstract}

\section{Introduction}

Partitions of subsets of $n$-dimensional Euclidean space into parts of smaller diameter are the subject of a large number of works related to the Borsuk conjecture~\cite{borsuk,Raigorodskii2007,zong2}. This conjecture states that every bounded set in $\R^n$ of positive diameter can be partitioned into $n+1$ parts of smaller diameter. The conjecture has been proved in a number of special cases; in particular, Hadwiger~\cite{Hadwiger1945} proved it for convex bodies with smooth boundary. In 1993, Kahn and Kalai proved that the conjecture is false for $n=1325$ and $n\geq 2015$~\cite{KahnKalai1993}. Recently, OpenAI announced\footnote{At present, there is no full formalization in Lean 4 or an independent review of this paper.} a counterexample in dimension~9~\cite{OpenAI}. Before that, the smallest dimension in which a counterexample was known was 64~\cite{bondarenko2014borsuk,jenrich201464}, and then 63, obtained by an AI-generated refinement of the Bondarenko's construction~\cite{ji2026borsuk63}. See~\cite{raigorodskii2008three} for a history of other improvements.

Quantitative versions of the problem have also been studied in a number of papers~\cite{Filimonov,kupavskii_d_3_5,DAM_2022_Tolmachev}. Let all sets under consideration have unit diameter, and fix $k \geq n+1$. If a partition into $k$ parts of smaller diameter is possible, how small can the largest diameter of the parts be made? We denote this number by $\beta_{n,k}$. Since the regular $n$-simplex cannot be partitioned into $n$ parts of smaller diameter, for $k=n$ we would have $\beta_{n,n}=1$; thus, in every dimension the first nontrivial case is $k=n+1$, the case of the Borsuk conjecture. It is known that $\beta_{1,2}=1/2$, $\beta_{2,3}=\sqrt{3}/2$~\cite{gale1953inscribing}, and
\[
0.88807\ldots = \sqrt{\frac{3+\sqrt{3}}{6}} \leq \beta_{3,4} \leq 0.9654\ldots
\]

The upper bound is obtained as follows. Makeev~\cite{makeev2000} showed that every set of unit diameter in $\R^3$ can be covered by a truncated rhombic dodecahedron with unit distance between each pair of parallel faces, and pointed out that this polyhedron can be partitioned into four parts of diameter $0.98$. The algorithm presented in~\cite{DAM_2022_Tolmachev} improves this bound somewhat, but with a single polyhedron it only yields $\beta_{3,4}\leq 0.9697\dots$. To obtain the bound $\beta_{3,4}\leq 0.966$ stated in~\cite{DAM_2022_Tolmachev}, it suffices to apply the same algorithm to two polyhedra obtained by a further truncation of Makeev's universal cover. A more careful optimization of the parameters gives $\beta_{3,4}\leq 0.9654$~\cite{borsuk_repo}, which is still far from the lower bound.

The main result of the present paper is the construction of a set that improves the lower bound for $\beta_{3,4}$, namely, a finite set $V\subset\R^3$ such that every partition of $V$ into four parts contains a part whose diameter is slightly larger than the corresponding value $0.88807\ldots\cdot\diam V$ for the sphere. We also present similar examples for $4 \leq n \leq 6$. All these examples were found by a global search starting from random initial configurations.

\paragraph{AI disclosure.} This work began in 2024; its original goal was to find a low-dimensional counterexample to the Borsuk conjecture. In April 2026, we carried out extensive code-evolution experiments with the GigaEvo framework~\cite{khrulkov2025gigaevo}\footnote{An open-source evolutionary framework for automating ML and LLM-oriented tasks developed in AIRI.} using DeepSeek~4, Qwen~3 Max, Gemini~3, and Claude Opus~4.6; the best ratio obtained in $\R^3$ was about $0.872$. In September 2026, we returned to the problem and, with the help of GPT-6 Astra, developed new algorithms and heuristics, some of which were proposed by the AI itself. GPT-6 was also used to prepare the GitHub repository~\cite{gale_repo} and to write the document describing the algorithms; Claude Opus 5.5 was used for proofreading and language editing.

\section{Preliminaries}\label{sec:prelim}

Throughout the paper, we consider only the Euclidean space $\R^n$. The Euclidean distance between points $x,y \in \R^n$ is denoted by $\|x-y\|$, and the diameter of a bounded set $A \subset \R^n$ by
\[
    \diam A = \sup_{x,y\in A} \|x-y\|.
\]

Fix a bounded set $A$ of positive diameter and partition it into $k$ disjoint subsets $A_1, \dots , A_k$, where $k\geq n+1$. We ask how small the largest of the diameters of the sets $A_i$ can be compared with the diameter of $A$.

\begin{definition}
    The \emph{diameter ratio} of $A$ for $k$ parts is
\[
        \beta_{k}(A) = \inf_{A_1 \sqcup \dots \sqcup A_k = A } \frac{\max_i \diam A_i}{\diam A} = \inf_{A_1 \cup \dots \cup A_k \supseteq A } \frac{\max_i \diam A_i}{\diam A}.
\]
\end{definition}
The second equality holds because a covering $A_1,\dots,A_k$ of $A$ yields the partition $A_i' = (A_i\cap A)\setminus(A_1\cup\dots\cup A_{i-1})$ of $A$ with $\diam A_i'\leq\diam A_i$.

\begin{definition}
For every pair $(n,k)$ with $k \geq n+1$, we define the extremal value
\[
    \beta_{n,k} = \max\limits_{A \subset \R^n, \; 0< \diam A<\infty} \beta_{k}(A).
\]
\end{definition}

A counterexample to the Borsuk conjecture in dimension $n$ exists if and only if $\beta_{n,n+1}=1$. Hence, the key open question can be stated as follows.

\paragraph{Minimal dimension of a Borsuk counterexample.} Determine the smallest $n$ such that $\beta_{n,n+1} = 1$.

The following problem, like many problems of geometric optimization, is unlikely to be solvable in full generality.

\paragraph{The quantitative Borsuk problem.} For given $n$ and $k$, determine $\beta_{n,k}$ and the class $\mathcal{A}_{n,k}$ of extremal closed convex sets on which this value is attained.

Traditionally, the case $k=n+1$ is of the greatest interest. It is known that $\beta_{2,3} = \sqrt{3}/2$~\cite{gale1953inscribing}.

\begin{theorem}[Lenz~\cite{Lenz}]
  $A \in \mathcal{A}_{2,3}$ if and only if $A$ is a disk.
\end{theorem}

This leads to the following conjecture, which is much less known than the Borsuk conjecture. Denote by $B^n$ the $n$-dimensional ball; note that $\beta_{n+1}(B^n)=\beta_{n+1}(S^{n-1})$.

\begin{conjecture}[D.~Gale, 1953~\cite{gale1953inscribing}]
    For every $n\geq 2$, the ball is an extremal set, i.e., $\beta_{n,n+1}=\beta_{n+1}(B^n)$.
\end{conjecture}

A stronger form of this conjecture, generalizing Lenz's theorem, asserts that $A \in \mathcal{A}_{n,n+1}$ if and only if $A$ is a ball.

Central projection of the facets of a regular $n$-simplex inscribed in $S^{n-1}$ partitions the sphere into $n+1$ congruent regular spherical simplices (taking cones from the center gives the corresponding partition of the ball). The diameter of each part is the distance between the centers of two complementary faces of dimensions $\lfloor n/2\rfloor-1$ and $\lceil n/2\rceil-1$: for even $n$, these are two opposite faces of dimension $n/2-1$; for odd $n$, faces of dimensions $(n-3)/2$ and $(n-1)/2$ (for $n=3$, a vertex of a spherical triangle and the midpoint of the opposite side). This gives the following upper bound for $n \geq 2$:
\[
\beta_{n+1}(B^n) = \beta_{n+1}(S^{n-1})\le \sigma_n :=
\begin{cases}
\displaystyle \sqrt{\frac{n+1}{n+2}},
& n \text{ even},\\[8pt]
\displaystyle
\sqrt{\frac{1}{2}+\frac{1}{2}\sqrt{\frac{n-1}{n+3}}},
& n \text{ odd}.
\end{cases}
\]
This partition is known to be optimal only for $n=2$ and $n=3$~\cite{kalai2015some}, i.e., $\beta_3(B^2)=\sigma_2=\sqrt{3}/2$ and $\beta_4(B^3)=\sigma_3=\sqrt{(3+\sqrt{3})/6}$.

Every counterexample to the Borsuk conjecture is also a counterexample to Gale's conjecture, but counterexamples to Gale's conjecture had not previously been constructed independently of the Borsuk problem; in particular, none were known in dimensions below~9. In this paper, we present a counterexample to Gale's conjecture in dimension~3, found numerically and verified in exact integer arithmetic.

\begin{theorem}\label{thm:main}
    There exists a set $V \subset \R^3$ with $|V|=101$ such that $\beta_4(V) > 0.88928 > \beta_4(S^2)$.
\end{theorem}

The first inequality can be verified in exact arithmetic. The second one is the bound $\beta_4(S^2)\leq\sigma_3=0.88807\ldots$ given by the partition of $S^2$ into four congruent spherical triangles obtained by central projection of the faces of an inscribed regular tetrahedron.

\section{Algorithms}

\subsection{Basic algorithm}

The starting point of our search for finite sets $V$ with large values of $\beta_k(V)$ was the following procedure.

\begin{enumerate}
\item Generate a random set $V = \{v_1, \dots, v_s\} \subset [0,1)^n$, where $s > k+1$.

\item Sort the pairs of indices $(i,j)$, $i<j$, in decreasing order of the distance $\|v_i - v_j\|$.

\item For $0\leq t\leq1$, let $G_t=([s],E_t)$ be the graph in which $(i, j)$ is an edge if and only if $\|v_i-v_j\|>t \diam V$. For $t=1$, the graph $G_t$ is empty, and for $t=0$, it is complete. As $t$ decreases, edges are added to $G_t$ in the order of the list sorted in step~2, so $\chi(G_t)$ is non-increasing in $t$. Using binary search over this list, we find its shortest initial segment $E$ such that the graph $([s],E)$ has chromatic number at least $k+1$; let $t^*\diam V$ be the length of the last edge in $E$. At each step of the binary search, the $k$-colorability of the current graph is decided by a SAT solver. The edge set $E$ is then fixed.

\item Using a finite-dimensional local optimization method, find a local maximum $q^*$ of the problem
\begin{equation}\label{eq:objective}
    \frac{\min\limits_{(i,j) \in E} \|v_i-v_j\|}{\diam V} \longrightarrow \max
\end{equation}
over the positions of the points $v_1,\dots,v_s$, with the edge set $E$ fixed.
\end{enumerate}

Since the graph $([s],E)$ is not $k$-colorable, every partition of the optimized set $V$ into $k$ parts contains a part with both endpoints of some edge of $E$. Hence $\beta_k(V)\geq q^*$, and the algorithm yields the bound $\beta_{n,k} \geq q^*$.

\subsection{Improvements}

We briefly list the most important modifications that eventually allowed us to find the desired finite set.

\paragraph{Reduction to the $k$-core.} We repeatedly delete the edges incident to vertices of degree less than $k$ (such vertices become isolated) until every non-isolated vertex has degree at least $k$. This does not affect the property $\chi\geq k+1$, since a vertex of degree less than $k$ can always be colored after its neighbors.

\paragraph{Distance graph pruning.} We extend the list of candidate edges using the condition $\ell_{ij}=\|v_i-v_j\| > (1-\delta)\,t^*\diam V$, where $\delta>0$ is a parameter of the algorithm. Then we add random noise: $\ell'_{ij}= \ell_{ij}+\eta Z_{ij}$, where $\eta>0$ is the noise amplitude and $Z_{ij}$ are independent random variables with the standard normal distribution $\mathcal{N}(0,1)$. We sort the edges by the perturbed lengths $\ell'_{ij}$ and, as in step~3, find an initial segment $E'$ of this list such that the graph $G'=([s],E')$ has chromatic number at least $k+1$.

\paragraph{Vertex splitting.} We add to the graph $G$ a twin $v'$ of some vertex $v$, i.e., a new vertex adjacent to the same vertices as $v$ and located at the same point. Then we delete some edges incident to $v$ and some edges incident to $v'$ in such a way that the resulting graph $G''$ still has chromatic number at least $k+1$. Next, we perturb the positions of $v$ and $v'$ and apply local optimization. If the value of the objective function~\eqref{eq:objective} improves, the split is accepted.

\paragraph{Paired splits.} Two vertices are split simultaneously, subject to the same conditions. In this case, the search becomes considerably more expensive.

\paragraph{Geometric refinement.} Using an empirical hypothesis about the piecewise smooth boundary of the convex set on which the points lie, we increase the number of points in $V$ or resample this boundary in a different way.

\medskip
As the main nonlinear optimization methods, we use L-BFGS-B, the limited-memory Broyden--Fletcher--Goldfarb--Shanno algorithm with bound constraints, and Sequential Least Squares Programming (SLSQP), both as implemented in SciPy. A detailed description of these and other algorithms used is given in the AI-written document in the repository~\cite{gale_repo}.

\section{Results}\label{sec:results}

With the basic algorithm, for $n=3$, $k=4$, and $30\leq s \leq 70$, we were able to generate configurations with $q^*= 0.865 \dots$. After the code was evolved with evolutionary algorithms, some of the changes (in particular, random perturbations of the distance graph) made it possible to reach $q^*= 0.872 \dots$. On the one hand, it is easy to obtain a finite set with a larger value of the objective function, since the sphere can be approximated by a finite set of points.

\begin{definition}
    A set $N\subseteq X$ is called an \emph{$\varepsilon$-net} in a metric space $M=(X,\rho)$ if for every $x \in X$ there exists $y \in N$ such that $\rho(x,y)\leq \varepsilon$.
\end{definition}

\begin{claim}
    Let $V$ be an $\varepsilon$-net on the sphere $S^2 \subset \R^3$ of unit diameter. Then $\beta_4(V) \geq \beta_4(S^2) - 2\varepsilon$. If, moreover, $\diam V=1$ (for instance, if $V$ is centrally symmetric), then $\beta_4(V) \leq \beta_4(S^2)$.
\end{claim}

On the other hand, the sets we obtained contained far fewer points than an $\varepsilon$-net giving the same bound. This gave hope that these methods could produce a set with $\beta_4(V) > \beta_4(S^2)$.

After the release of GPT-6 Astra, several auxiliary algorithms were developed within a few days, and they indeed made it possible to find such sets. Some properties of the constructed sets may be surprising:
\begin{itemize}
\item The value $\beta_4(V)$ exceeds $\beta_4(S^2)$ only by about $0.001$. If further improvements are possible, they require much greater computational resources or new ideas.

\item Attempts to replace the finite set $V$ by a body $F$ with piecewise-spherical of parabolic boundary lead to a smaller value of $\beta_4(F)$.

\item The best sets in dimension~$3$ with $50$ to $250$ points, selected from several thousand restarts, have a similar structure.
\end{itemize}
%At present, we cannot offer an explanation of these phenomena.

The code and the configurations found are available in the repository~\cite{gale_repo}. To verify the 101-point configuration (Fig.~\ref{fig:101p}; the coordinates are given in Appendix~\ref{app:coords}), we short independent verivication scripts: it suffices to compute the diameter $D$ of the point set, select the pairs of points at distance greater than $0.88928\,D$, and check that the corresponding distance graph admits no proper $4$-coloring.

The best configurations obtained in independent runs all have the following structure: three long arcs that divide the equator into unequal parts, two triples of short arcs near the poles, and several points inside the trihedral angles at each pole. Imposing exact symmetry, for example requiring the long arcs to be arranged symmetrically, slightly worsens the result.

\begin{figure}[ht]
    \centering
    \includegraphics[width=0.45\linewidth]{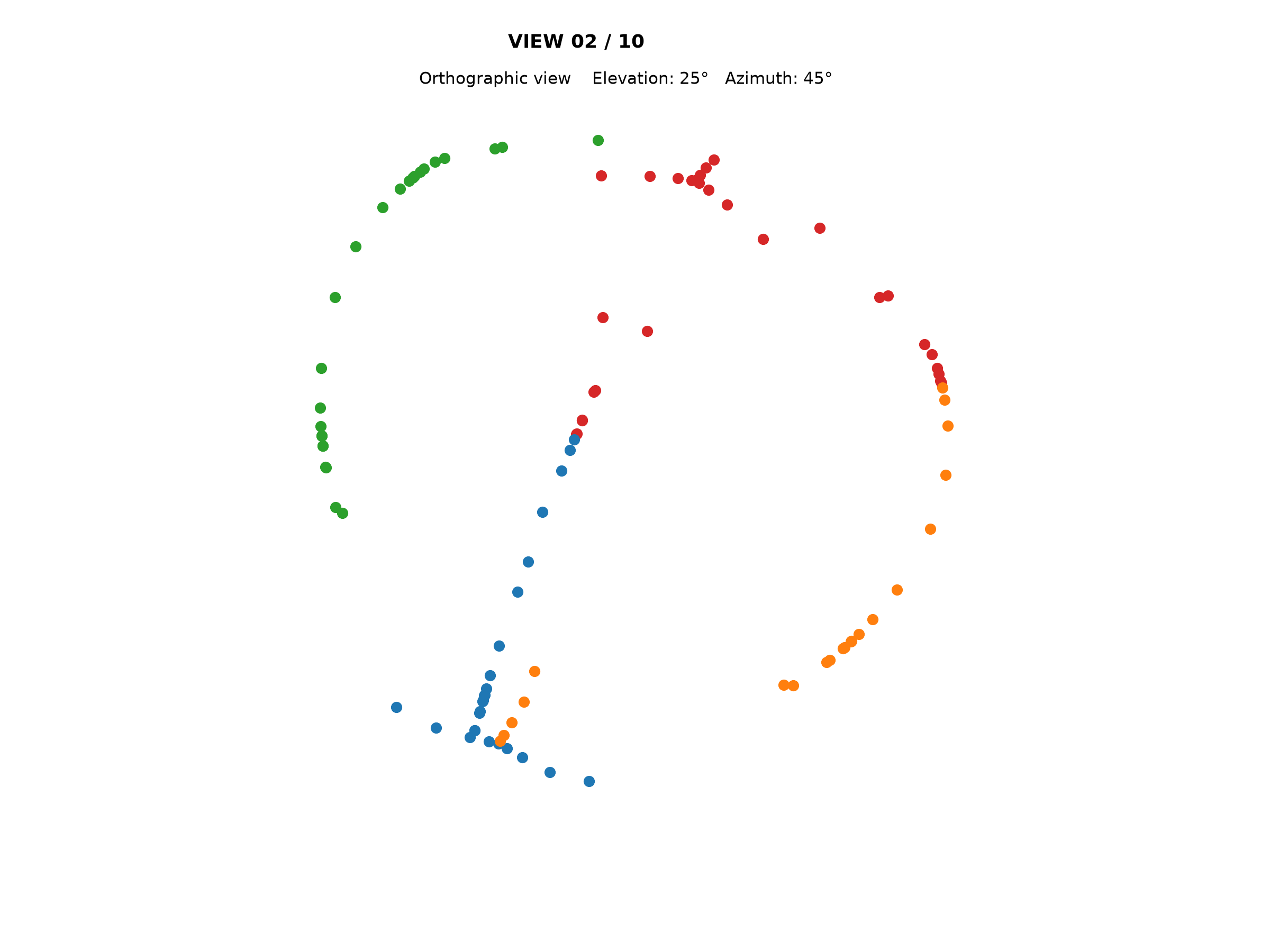}\hfill
    \includegraphics[width=0.45\linewidth]{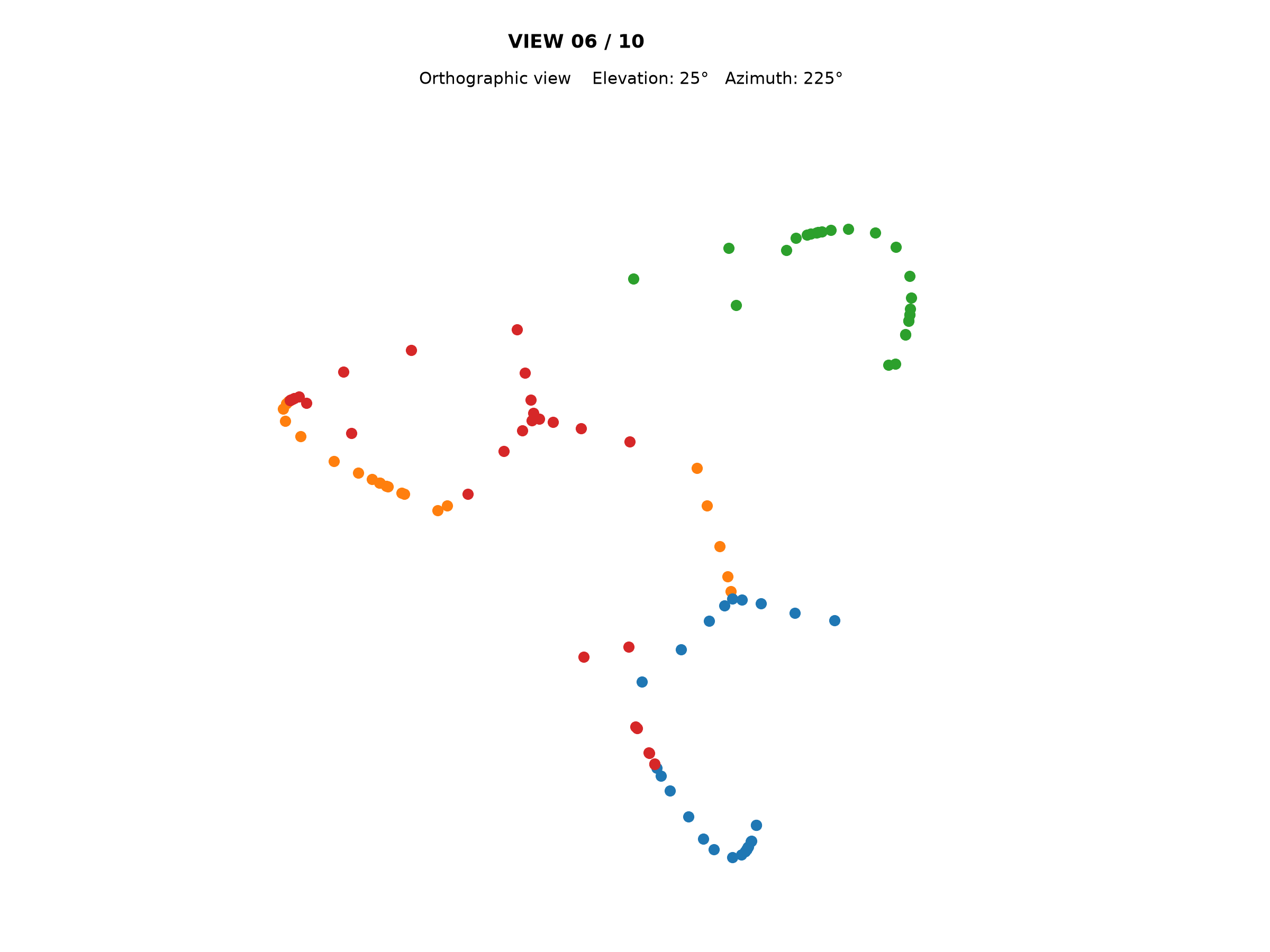}\\[1ex]
    \includegraphics[width=0.45\linewidth]{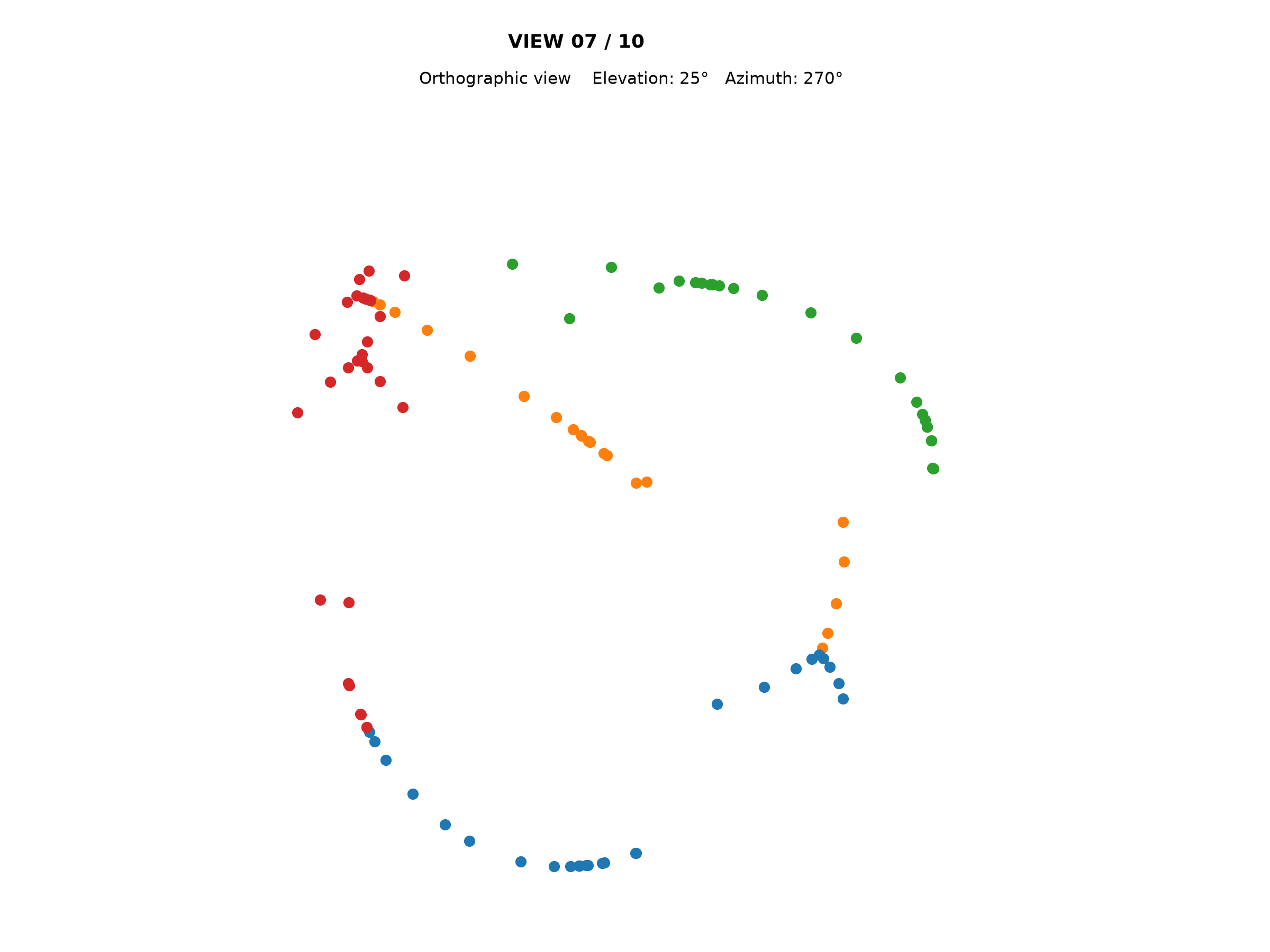}\hfill
    \includegraphics[width=0.45\linewidth]{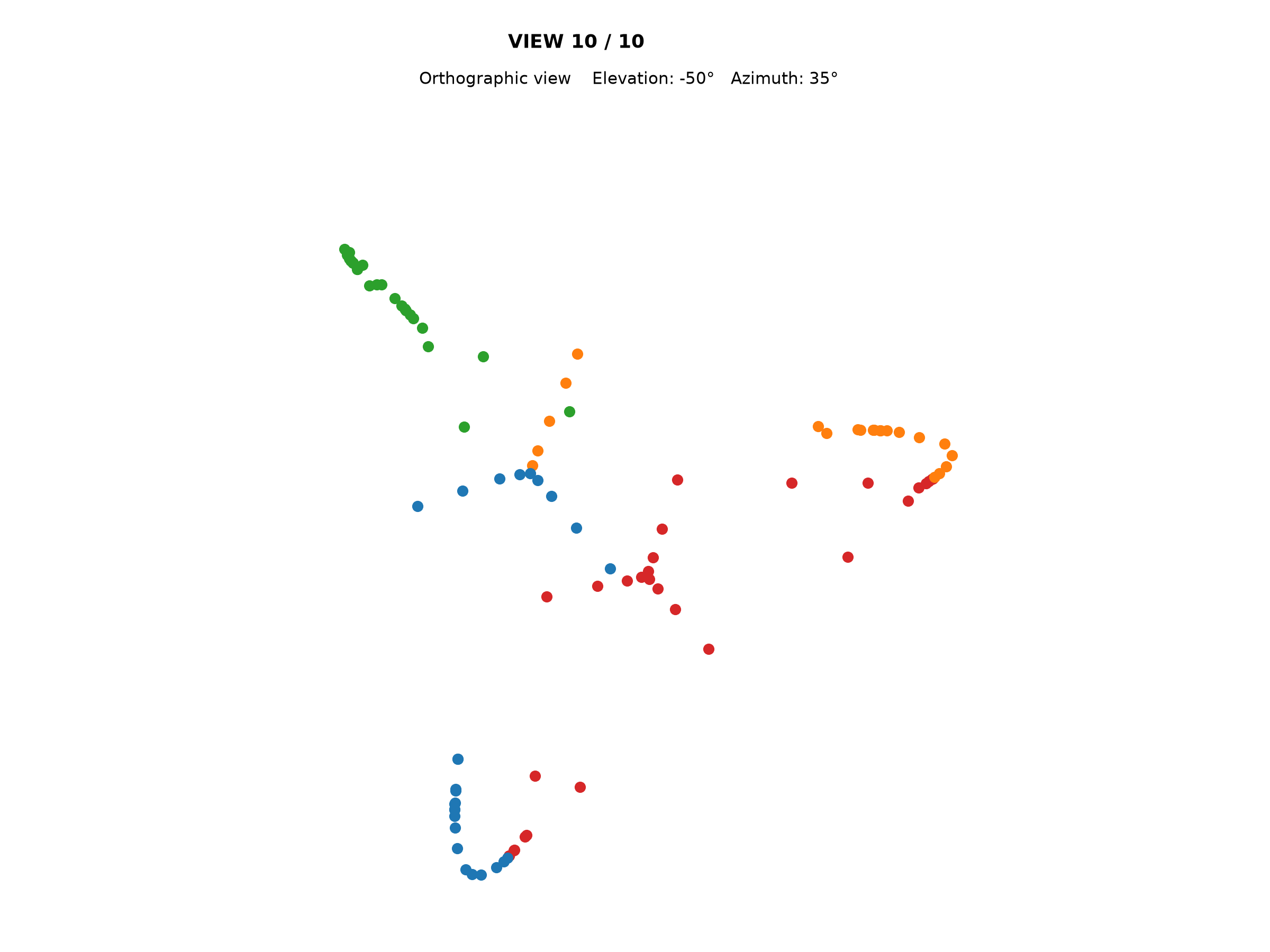}
    \caption{The 101-point configuration in Theorem~\ref{thm:main}.}
    \label{fig:101p}
\end{figure}

The repository also contains counterexamples to Gale's conjecture in dimensions 4, 5, and 6 (see Table~\ref{tab:results}). These values are most likely far from optimal, but in each case they exceed the upper bound $\sigma_n$ for the ball.

\begin{table}[ht]
    \centering
    \begin{tabular}{|c|c|c|c|}
        \hline
        $n$ & $\sigma_n$ & $\beta_{n+1}(V)$ & $|V|$ \\
        \hline
        3 & $0.888073\ldots$ & $0.88928\ldots$ & 101 \\
        4 & $0.912870\dots$ & $0.91314\ldots$ & 61 \\
        5 & $0.923879\dots$ & $0.92576\ldots$ & 87 \\
        6 & $0.935414\dots$ & $0.93696\ldots$ & 64 \\
%        7 & $0.941965\dots$ & & \\
        \hline
    \end{tabular}
    \caption{Current records: the diameter ratio $\beta_{n+1}(V)$ of the best sets $V\subset\R^n$ found, compared with the upper bound $\sigma_n\geq\beta_{n+1}(B^n)$ given by the simplex partition.}
    \label{tab:results}
\end{table}

\section{Discussion}

In G.~Kalai's blog~\cite{kalai_blog}, as well as in the survey~\cite{kalai2015some}, the following question was raised: does the smallest dimension in which a counterexample to the Borsuk conjecture exists depend on whether the counterexample is required to be a finite set? This question remains open. The smallest dimension in which a finite counterexample is known is still 63. It can be shown that every finite set of points chosen on the embedding of $\mathbb{RP}^3$ into $\R^9$ constructed in~\cite{OpenAI} can be partitioned into 8 parts of smaller diameter.

This question can be generalized as follows.

\begin{question}
    For which values of the parameters $n,k$ is the value $\beta_{n,k}$ attained on a finite set?
\end{question}

One can also regard the number of points $m$ as a parameter and consider
\[
\beta_{n,k,m}=\max\{\beta_k(V) : V\subset\R^n,\ |V|=m\},
\]
which gives a generalization of the min-max ratio problem~\cite{bateman1951geometrical,berthold2026global}. The following question is beyond the scope of this paper.

\begin{question}
    What is the smallest $m=m(n)$ for which there exists an $m$-point set that is a counterexample to Gale's conjecture in $\R^n$?
\end{question}

At present, it is unclear whether there exist counterexamples whose properties can be established by more conceptual mathematical arguments than an exhaustive computer search or a SAT solver.

\begin{problem}
    Prove that a counterexample to Gale's conjecture exists for every $n \geq 3$. Describe the construction geometrically and prove the statement without exhaustive enumeration of graph colorings.
\end{problem}

The examples presented in this paper are substantially non-spherical: if the points are radially projected onto a sphere centered at the center of the smallest enclosing ball, the resulting set is no longer a counterexample. Note that a counterexample $V\subset S^{n-1}$ cannot contain a pair of antipodal points, since otherwise $\diam V=\diam S^{n-1}$ and $\beta_{n+1}(V)\leq\beta_{n+1}(S^{n-1})$. On the other hand, the counterexamples to the Borsuk conjecture constructed in~\cite{KahnKalai1993,bondarenko2014borsuk} consist of vectors of equal length, so spherical counterexamples to Gale's conjecture exist in sufficiently high dimensions.

\begin{question}[see~\cite{kalai2015some}]
What is the smallest $n$ for which there exists a counterexample to Gale's conjecture that is a subset of $S^{n-1}$?
\end{question}

The following question is a long-standing open problem.

\begin{question}[see~\cite{kalai2015some}]
    Is the partition of $S^{n-1}$ into $n+1$ regular spherical simplices optimal for every $n\geq 4$, i.e., is $\beta_{n+1}(B^n)=\sigma_n$?
\end{question}

Of course, the central question remains whether there exists a counterexample to the Borsuk conjecture in dimension lower than $9$. It seems that the ideas and algorithms described in this paper are not sufficient to find such counterexamples by multistart global optimization. On the other hand, the bound $\beta_{4,8} \leq 0.9999$ obtained in~\cite{tolmachev2026reducing} indirectly suggests that the Borsuk conjecture in dimension $4$ would be very difficult to prove computationally, which leaves open the possibility that a counterexample exists already for $n=4$.

\paragraph{Acknowledgements.} The author is grateful to A.~M.~Raigorodskii for drawing attention to this problem.

\appendix
\section{Coordinates of the 101-point counterexample}\label{app:coords}

Table~\ref{tab:borsuk-101-coordinates} lists the points of the set $V$ from Theorem~\ref{thm:main}. Theorem~\ref{thm:main} holds for these rounded coordinates themselves: if $D$ is the diameter of the set, then the $2030$ pairs of points at distance greater than $0.88928\,D$ form a graph that admits no proper $4$-coloring.

% Canonical source: SELECTED_3d.zip / 3d_088891/points.csv
% SHA-256: 362e063b96d5075c6904446d3e9164879b47233de32158e46ff6703c96704093
% Original point order, read left to right and then top to bottom.
% All coordinates rounded to seven digits AFTER the decimal point.
% No additional LaTeX packages are required by this table.
\begin{table}[!ht]
\centering
\footnotesize
\setlength{\tabcolsep}{5pt}
\renewcommand{\arraystretch}{1.0}
\caption{101-point configuration in $\R^3$.}
\label{tab:borsuk-101-coordinates}
\begin{tabular}{rrr||rrr}
\hline
$x_1$ & $x_2$ & $x_3$ & $x_1$ & $x_2$ & $x_3$ \\
\hline
$0.3807914$ & $0.0962028$ & $-0.3242537$ & $-0.3560793$ & $-0.1353715$ & $0.3102441$ \\
$-0.2852755$ & $-0.2804681$ & $0.2982384$ & $-0.4495227$ & $-0.0877076$ & $0.1986231$ \\
$-0.0235972$ & $-0.2726832$ & $-0.4958457$ & $0.3768454$ & $0.1660764$ & $-0.2992439$ \\
$0.5407728$ & $-0.0252179$ & $0.0744869$ & $0.3941294$ & $0.0350691$ & $-0.3240403$ \\
$-0.2191007$ & $-0.3760150$ & $-0.3755289$ & $-0.3695321$ & $-0.3618959$ & $-0.0000334$ \\
$-0.0199093$ & $0.5528233$ & $-0.1283863$ & $-0.0491757$ & $-0.2900877$ & $-0.4874451$ \\
$-0.4221627$ & $0.1966422$ & $0.2009433$ & $0.3686276$ & $0.1503312$ & $-0.3167916$ \\
$-0.3292606$ & $-0.3931751$ & $-0.2252822$ & $0.5154236$ & $-0.1001790$ & $0.2226708$ \\
$-0.0253103$ & $-0.2554733$ & $0.4393427$ & $-0.1803981$ & $0.5503760$ & $-0.0009600$ \\
$-0.0099454$ & $-0.2630252$ & $-0.4995430$ & $0.4007449$ & $0.2578140$ & $-0.1504765$ \\
$0.3505919$ & $-0.1865142$ & $0.4174207$ & $0.2050635$ & $0.1831177$ & $-0.4293148$ \\
$0.2746356$ & $-0.2025288$ & $0.4546163$ & $0.3640042$ & $0.1422828$ & $-0.3256709$ \\
$-0.3211546$ & $0.4411104$ & $0.1381561$ & $0.4216235$ & $-0.1610749$ & $0.3615418$ \\
$-0.1816098$ & $-0.3620061$ & $-0.4102977$ & $-0.3373248$ & $-0.3919925$ & $-0.2094446$ \\
$-0.3488780$ & $-0.1319739$ & $0.3197610$ & $-0.3987962$ & $-0.1159669$ & $0.2650061$ \\
$-0.0458995$ & $0.5578411$ & $-0.1102685$ & $0.0784715$ & $0.5070985$ & $-0.1991818$ \\
$0.2784339$ & $0.1700810$ & $-0.3940311$ & $0.5247249$ & $-0.0890238$ & $0.1968562$ \\
$0.4030776$ & $0.2373564$ & $-0.2097259$ & $-0.3703436$ & $-0.1291660$ & $0.2961362$ \\
$0.2075988$ & $-0.2113823$ & $0.4751283$ & $0.4008200$ & $-0.0461820$ & $-0.3126089$ \\
$-0.2476206$ & $0.5170097$ & $0.0589887$ & $0.4904031$ & $-0.1221274$ & $0.2753791$ \\
$-0.2977037$ & $0.4711859$ & $0.1114675$ & $-0.3404125$ & $-0.1171323$ & $0.3348009$ \\
$0.0398331$ & $-0.1740864$ & $0.4899806$ & $-0.3116044$ & $-0.3951525$ & $-0.2558847$ \\
$-0.2833496$ & $-0.0304706$ & $0.4085761$ & $-0.3404423$ & $-0.1667078$ & $0.3129268$ \\
$0.3523596$ & $0.1485984$ & $-0.3356124$ & $-0.4138401$ & $-0.3074944$ & $-0.0205616$ \\
$-0.1149363$ & $-0.1178070$ & $0.4690417$ & $0.1138774$ & $-0.2265340$ & $0.4783584$ \\
$-0.3380685$ & $0.1489373$ & $0.3323502$ & $-0.3575756$ & $0.3766951$ & $0.1833862$ \\
$0.5326245$ & $-0.0770263$ & $0.1696352$ & $0.0283772$ & $0.5365805$ & $-0.1622242$ \\
$0.2299305$ & $-0.2089227$ & $0.4693788$ & $0.3276424$ & $0.1579222$ & $-0.3566089$ \\
$-0.3689161$ & $-0.3804383$ & $-0.1344716$ & $-0.3476660$ & $0.3974084$ & $0.1701033$ \\
$-0.3376319$ & $0.4154322$ & $0.1578776$ & $0.3705793$ & $0.1275514$ & $-0.3244517$ \\
$0.1712134$ & $-0.2146969$ & $0.4821205$ & $-0.3315974$ & $0.4253427$ & $0.1506028$ \\
$0.1984536$ & $-0.2122531$ & $0.4771821$ & $0.0779847$ & $-0.1959417$ & $-0.5090523$ \\
$-0.3411717$ & $-0.3912094$ & $-0.2014937$ & $0.0032233$ & $-0.2535311$ & $-0.5026568$ \\
$-0.0110774$ & $-0.2638209$ & $-0.4992591$ & $0.5286694$ & $-0.0834042$ & $0.1840350$ \\
$-0.3512779$ & $-0.3883896$ & $-0.1793195$ & $-0.0061716$ & $0.5488859$ & $-0.1382818$ \\
$-0.3489005$ & $-0.1448964$ & $0.3140760$ & $-0.0072590$ & $0.5492202$ & $-0.1375054$ \\
$0.5391917$ & $-0.0636072$ & $0.1393137$ & $0.3907369$ & $0.1975685$ & $-0.2624424$ \\
$-0.2701644$ & $-0.3952287$ & $-0.3136478$ & $-0.3211385$ & $-0.0858187$ & $0.3642750$ \\
$-0.0966865$ & $0.5602195$ & $-0.0748505$ & $-0.3534674$ & $0.2830232$ & $0.2554356$ \\
$-0.3213532$ & $-0.2095827$ & $0.3095328$ & $-0.1012565$ & $-0.3225157$ & $-0.4642462$ \\
$0.0290344$ & $-0.2342847$ & $-0.5074368$ & $0.0332040$ & $0.5346343$ & $-0.1651591$ \\
$0.1455650$ & $-0.2161549$ & $0.4855848$ & $0.0949415$ & $0.5019801$ & $-0.1949247$ \\
$0.0064962$ & $0.5448308$ & $-0.1471536$ & $0.0257652$ & $-0.2366885$ & $-0.5070009$ \\
$0.0007126$ & $-0.2553278$ & $-0.5021253$ & $0.0043313$ & $0.5455470$ & $-0.1456805$ \\
$0.5422049$ & $-0.0386625$ & $0.0800017$ & $-0.3717280$ & $0.3462640$ & $0.1865461$ \\
$-0.3429776$ & $0.4059802$ & $0.1645686$ & $0.1809769$ & $-0.2137626$ & $0.4806143$ \\
$0.1942764$ & $-0.2126746$ & $0.4780348$ & $-0.3354618$ & $0.4190416$ & $0.1553035$ \\
$-0.3703952$ & $-0.3785746$ & $-0.1316978$ & $-0.3507379$ & $-0.3885884$ & $-0.1805062$ \\
$-0.3414434$ & $-0.3911424$ & $-0.2009338$ & $0.5325432$ & $-0.0771622$ & $0.1699623$ \\
$0.5390925$ & $-0.0638373$ & $0.1399082$ & $0.0780330$ & $-0.1958849$ & $-0.5090517$ \\
$0.5286521$ & $-0.0834295$ & $0.1840944$ &  &  &  \\
\hline
\end{tabular}
\end{table}

\clearpage
\bibliographystyle{abbrv}
\bibliography{main1}

\end{document}